\documentclass[pre,twocolumn,twoside,showpacs,superscriptaddress,floatfix]{revtex4-1}

\usepackage{graphics,amssymb,amsmath,epsf,color,hyperref}
\epsfclipon
\usepackage{epsfig,bm,epsf,graphics}
\usepackage{graphicx}
\usepackage{dcolumn}
\usepackage{bm}
\usepackage{braket}
\usepackage{color}
\usepackage{soul}
\usepackage{floatrow}
\usepackage{hyperref}
\usepackage{multirow}

\hypersetup{colorlinks, citecolor=blue, linkcolor=blue}

\usepackage[table]{xcolor}
\usepackage{amsmath}
\usepackage{subfigure}

\newcommand{\dd}{\mathrm{d}}
\renewcommand{\root}{\mathrm{root}}
\newcommand{\singular}{\mathrm{s}}
\newcommand{\x}{\langle x \rangle}
\begin{document}

\title{Annealing approach to root-finding}

\author{Junghyo Jo}
\email{jojunghyo@snu.ac.kr}
\affiliation{Department of Physics Education, Seoul National University, Seoul 08826, Korea}
\affiliation{Center for Theoretical Physics and Artificial Intelligence Institute, Seoul National University, Seoul 08826, Korea}
\affiliation{School of Computational Sciences, Korea Institute for Advanced Study, Seoul 02455, Korea}

\author{Alexandre Wagemakers}
\email{alexandre.wagemakers@urjc.es}
\affiliation{Nonlinear Dynamics, Chaos and Complex Systems Group, Departamento de Física, Universidad Rey Juan Carlos, Tulipán s/n, Móstoles 28933, Madrid, Spain}

\author{Vipul Periwal}
\email{vipulp@mail.nih.gov}
\affiliation{Laboratory of Biological Modeling, National Institute of Diabetes and Digestive and Kidney Diseases, National Institutes of Health, Bethesda, Maryland 20892, USA}

\date{\today}

\begin{abstract}
The Newton-Raphson method is a fundamental root-finding technique with numerous applications in physics.
In this study, we propose a parameterized variant of the Newton-Raphson method, inspired by principles from physics. 
Through analytical and empirical validation, we demonstrate that this novel approach offers increased robustness and faster convergence during root-finding iterations. Furthermore, we establish connections to the Adomian series method and provide a natural interpretation within a series framework. Remarkably, the introduced parameter, akin to a temperature variable, enables an annealing approach. This advancement sets the stage for a fresh exploration of numerical iterative root-finding methodologies.
\end{abstract}

 
\maketitle


\section{Introduction}
\label{intro}

Root-finding algorithms are important to solve equations, a fundamental task in quantitative theoretical science. Root-finding and fixed-point iterations are intricately connected, serving as essential tools in various fields such as optimization and algorithm development. 
In physics, root-finding methods can be smartly employed to find periodic orbits~\citep{abad2011computing, wang2018accelerated} and estimate parameters~\citep{amritkar2009estimating} in nonlinear systems with chaotic dynamics. 
Additionally, they can be used to identify stationary states of potential energy functions in classical systems~\citep{wales2003stationary, mehta2012energy}, in a generic Turing reaction-diffusion system~\citep{woolley2010analysis}, and in a nonlinear Schr\"odinger lattice~\citep{naether2011peierls}, as well as to find complex saddle points in quantum many-body systems~\citep{tomsovic2018complex}.
This application can be naturally extended to the maximization of log-likelihood~\citep{pfefferle2020whittle} and the minimization of free energy~\citep{maragakis2006optimal} of physical systems. 
They are also applicable in solving differential equations of fluid dynamics~\citep{herrada2016numerical} and self-consistent equations in quantum heat transport problems~\citep{saaskilahti2013thermal}.
Recently, the rise of deep learning has introduced complex, high-dimensional root-finding challenges, particularly in optimizing architecturally intricate neural networks, driving the need for further advancements in root-finding techniques.

Here, we introduce a concrete example how the root-finding method can be applied to solve a physics problem.
The Kuramoto model~\citep{kuramoto1975self} is paradigmatic in the theoretical study of the spontaneous synchronization of physical constituents. Consider a system of $N$ physical rotors with configurations $\theta_i,$ interacting as described by  
\begin{equation}\label{eq:Kuramoto}
    \frac{d\theta_i}{dt} = \kappa\sum_{j\not= i} \Gamma_{ij} \sin\left(\theta_j - \theta_i + \Psi_{ij}\right),
\end{equation}
where $\kappa$ is an overall coupling constant, $\Gamma_{ij}$ is a real weight matrix describing the edges of a weighted asymmetric directed graph with no multi-edges, and $\Psi_{ij}$ is a matrix of phase delays associated with the edge $j \rightarrow i.$ 
Note that we consider the case where no rotor has intrinsic angular velocity. 
By including slow learning dynamics for the $\Gamma$ and $\Psi$ matrices one could in principle store multiple patterns in the dynamics of this system, which may be relevant for modeling neuronal systems~\citep{buzsaki2004neuronal,kim2020dynamics}.
Now, suppose we consider the synchronization state of the rotors in a common angular velocity $\omega$, such that 
\begin{equation}\label{eq:ansatz}
    \theta_i(t) = \omega t  - \phi_i,
\end{equation}
where $\phi_i$ is a rotor-specific phase shift. Given $\kappa, \Gamma$ and $\Psi,$ is there a vector of phase shifts ${\phi}$ such that synchronization can occur for $\omega \neq 0?$ To exhibit this as a problem of finding roots for a set of $N-1$ coupled equations, inserting Eq.~(\ref{eq:ansatz}) into  Eq.~(\ref{eq:Kuramoto}) gives
\begin{equation}\label{eq:Kuramoto1}
    \omega = \kappa\sum_{j\not= i} \Gamma_{ij} \sin\left(\phi_i - \phi_j + \Psi_{ij}\right) \equiv \kappa f_i(\phi), \phantom{M}\forall i.
\end{equation}
We can set the initial phase shift \(\phi_0 = 0\) for the rotor with \(i = 0\).
Then $\omega(\phi)$ is specified by
\begin{equation}\label{eq:omega}
    \omega = \kappa\sum_{j\not= 0} \Gamma_{0j} \sin\left(-\phi_j + \Psi_{0j}\right),
\end{equation}
and inserting this value of $\omega$ in the remaining $N-1$ equations gives for all $N>i>0$:
\begin{equation}\label{eq:Kuramoto2}
    \sum_{j\not= 0} \Gamma_{0j} \sin\left(-\phi_j + \Psi_{0j}\right) = f_i(\phi).
\end{equation}
Now $\kappa$ does not appear in Eq.~(\ref{eq:Kuramoto2}) which is not surprising because changing $\kappa$ merely changes the scale of the time variable in this system. When the rotors have intrinsic individual frequencies, $\kappa$ cannot be scaled out and the Kuramoto synchronization phase transition occurs only for strong enough coupling.

To demonstrate that this problem is an exercise in root-finding, consider the case of two rotors. Then there is one equation determining $\phi_1:$
\begin{equation}
    \Gamma_{01}\sin(\phi_1 - \Psi_{01}) + \Gamma_{10} \sin(\phi_1 + \Psi_{10}) = 0.
\end{equation}
It is evident that for appropriate choices of coupling weights $\Gamma_{01},\Gamma_{10}$ and edge-specific phase shifts $\Psi_{01},\Psi_{10},$ this equation has a real root. In particular, the root then determines the synchronization frequency as shown in Eq.~(\ref{eq:omega}).

Given the numerous applications, root-finding methods historically predate even calculus, and continue to evolve, with well-known methods like bisection, regula falsi, Newton-Raphson, and secant methods~\citep{press2007numerical}. Notably, the Newton-Raphson method stands out as a powerful iterative approach for root-finding~\citep{ypma1995historical}.
The following iterative update of $x_n$,
\begin{equation}
    x_{n+1} = x_n - \frac{f(x_n)}{f'(x_n)},
\end{equation}
ultimately converges to the root, $x_\root$, satisfying $f(x_\root)=0$.
Note that $f'(x_n)$ represents the slope of $f(x)$ at $x=x_n$.

Research efforts persist in exploring hybrid methods, parallelized techniques, and symbolic and numerical hybrids, all aimed at enhancing efficiency, convergence, robustness, and global optimization. Physical analogies have long been important in numerical analysis. Can root-finding also be phrased in a physically motivated setting?
In this study, we develop an extended Newton-Raphson method, inspired by principles from physics. The method can be summarized in two iterative steps:
\begin{align}
\label{eq:newnewton}
    \hat{x}_{n+1} &= x_n - \frac{f(x_n)}{f'(x_n)}, \nonumber \\
    x_{n+1} &= \hat{x}_{n+1} - \beta \frac{f(\hat{x}_{n+1})}{f'(x_n)}.
\end{align}
Here, we introduce an auxiliary parameter $\beta$, where $\beta=0$ corresponds to the original Newton-Raphson method.
Although we present the method for a scalar variable, it can be readily extended to encompass vector functions with multiple variables.
For the specific application to the synchronization problem discussed above, we define
\begin{equation}\label{eq:Adef}
    A_{ik}^{-1}({\phi}) \equiv \frac{\partial}{\partial\phi_k} f_i({\phi}) .
\end{equation}
The iteration to find $\phi_i$ involves two steps:
\begin{equation}
    \hat\phi_k^{(n+1)} = \phi_k^{(n)} - \sum_j A_{kj}({\phi}^{(n)}) f_j({\phi}^{(n)})
\end{equation}
followed by
\begin{equation}
    \phi_k^{(n+1)} = \hat\phi_k^{(n+1)} - \beta\sum_j A_{kj}({\phi}^{(n)})  f_j({\hat\phi}^{(n+1)}).
\end{equation}
In fact, in this system there could well be multiple independent synchronization frequencies for subsets of rotors, depending on the structure of $\Gamma$ and $\Psi.$

This paper is organized as follows: In Section \ref{physical}, we provide a physical perspective on the Newton-Raphson method, highlighting its limitations and proposing solutions. Section \ref{sec:iterative} delves into the development of a parameterized variant of the Newton-Raphson method, elucidating its iterative implementation of Eq.~(\ref{eq:newnewton}) and demonstrating superior convergence and robustness through practical examples. Section \ref{sec:adomian} explores the connection between the extended Newton-Raphson method and the Adomian decomposition method, interpreting the auxiliary parameter $\beta$ as a natural expansion parameter for the Adomian series.
Finally, adopting the annealing concept in physics~\citep{kirkpatrick1983optimization}, we show further computational improvement by adjusting the $\beta$ value during iterations.

\section{Physical approach to finding roots}
\label{physical}
We now rederive the Newton-Raphson method from a novel perspective from physics, naturally extending it into an enhanced version.

\subsection{New perspective of the Newton-Raphson method}
Suppose that we have a scalar function $f(x)$ with a single variable $x$, aiming to determine a root, $x_\root$, such that $f(x_\root)=0$.
Then, we consider an integral, 
\begin{equation}
    Z \equiv \int \dd x \exp\left[-\frac{1}{2g^2} f(x)^2\right],
\end{equation}
where $g$ is a real-valued scale parameter.
We provide the mathematical justification for this form based on cohomological quantum field theory in 
Appendix \ref{topological}. For any function $h$ of $x,$ we expect that in the limit $g\downarrow 0,$
\begin{equation}
    \langle h \rangle \equiv \frac{1}{Z} \int \dd x \ h(x) \exp\left[-\frac{1}{2g^2} f(x)^2\right] \rightarrow h(x_\root),
    \label{eq:z}
\end{equation}
where $x_\root$ is a minimum or a root of $f,$ assuming a unique root or a unique minimum. In particular, the expectation value of $h(x)=x$ itself should be an estimate of $x_\root.$

The systematic approach to analyzing Eq.~(\ref{eq:z}) in the limit $g\downarrow 0$ is by using Laplace's approximation for the integral. In brief, this approximation requires us to find $x^*$ such that $f(x^*)f'(x^*) = 0,$ and then do a Taylor expansion about $x^*,$ evaluating the integral as an integral over fluctuations about $x^*.$ 

For the purposes of finding roots of $f,$ this is pointless as we have no idea where the root might be. Attempting to localize the integral at an arbitrary specific value $x_0$ using a Taylor expansion leads to
\begin{equation}
    Z = \int \dd x \exp\left[-\frac{1}{2g^2} \bigg(f(x_0) + f'(x_0)(x-x_0) + \ldots \bigg)^2\right].
\end{equation}
Defining a fluctuation variable {$\delta \equiv x-x_0$} with $f_0 \equiv f(x_0)$ and $f'_0 \equiv f'(x_0),$
\begin{equation}
    Z =  \int \dd\delta \exp\left[-\frac{1}{2g^2} \bigg(f_0 + f'_0\delta + \ldots \bigg)^2\right].
    \label{eq:Z}
\end{equation}
Then, we can find a better stationary point by varying
\begin{equation}
    \frac{1}{2g^2} \left[f_0 f'_0\delta + \frac{1}{2}(f'_0\delta)^2 \right],
    \label{eq:tovary}
\end{equation}
which is stationary when 
\begin{equation}
    \delta = -f_0/f'_0.
    \label{eq:eom}
\end{equation}
Since (up to higher order derivatives) this defines the expectation value of $\delta,$ it follows that the estimate for the root of $f$ given by this is
\begin{equation}
    x_\root \approx\x =  x_0 - f_0/f'_0,
\end{equation}
which is the Newton-Raphson update. 

As the basic assumption of the Laplace approximation is that the initial starting point is a stationary point of $f^2$ so that the expectation value of fluctuations vanishes, this attempt to expand the integral about an arbitrary value $x_0$ tells us that for a consistent Laplace approximation, we need to replace $x_0 \rightarrow x_0 -f_0/f'_0,$ and start over.

\subsection{Inconsistency in the Newton-Raphson method}

Laplace's approximation holds true only when we expand around a stationary point. As noted, the presence of a non-zero right-hand side in Eq.~(\ref{eq:eom}) indicates that this condition does not hold for arbitrary $x_0$.
However, there is another problem with our approach. We assumed that the linear approximation was actually solving the equation of interest, namely $f(x_\root) = 0$, but all we really solved was
\begin{equation}
    f_L(x) \equiv f_0 + f'_0 (x - x_0) = 0.
\end{equation}
This becomes apparent when we substitute the value $x_0 \rightarrow x_0 + \delta$ into Eq.~(\ref{eq:Z}), as the term $f_0^2/2g^2$ also undergoes modification, although it did not contribute to the variational equation, Eq.~(\ref{eq:tovary}).
Only if the corrections to Eq.~(\ref{eq:eom}) are of order $g$ and higher with $\delta=-f_0/f'_0=O(g)$ and $f_0^{''}\equiv f^{''}(x_0),$ 
\begin{equation}
    f(x_0 + \delta) \approx f_0 + f'_0 \delta + \frac{1}{2}  {{f^{''}_0}} \delta^2 = O(g^2),
\end{equation}
which is higher order in $g.$ 
This result implies that the Laplace approximation is na\"\i vely self-consistent. 
However, it is important to note that $f_0/f'_0=O(g) \Leftrightarrow f_0/g = O(f'_0).$ 
For an arbitrarily chosen point $x_0$, there is no assurance of such a relationship between $f_0$ and $f'_0$,
as $g$ must be as small as possible for the Laplace approximation to be valid. 
Hence, we deduce that the genuine expansion parameter is $f_0/g$, and Eq.~(\ref{eq:eom}) must be amended.

One approach to addressing both problems at the same time is to seek an enhanced integrand, often referred to as an effective action in analogous contexts within physics. Typically, such an effective integrand is computed by taking short distance scale fluctuations into account, leading to an integrand governing longer distance fluctuations but necessarily dependent on $g.$ For instance, \citet{lepage1993viability} demonstrated the reordering of strong-coupling terms to derive an improved effective action in lattice gauge theory. However, in our scenario, where only a single integral is involved, such a strategy cannot be directly applied. Nonetheless, we can explore the possibility of identifying a conceptually similar effective equation to refine the Laplace approximation, particularly in cases where the Taylor expansion employed to derive Eq.~(\ref{eq:eom}) becomes invalid.

$f_0/f'_0$ has the dimensions of length so fluctuations about $x_0$ are naturally considered big or small relative to this length scale. Two clues towards finding such an effective variant of Eq.~(\ref{eq:eom}) are the following:
\begin{itemize}
    \item When an exact stationary point for the Laplace approximation is not used, the leading correction to the logarithm of the integral is singular in a power of $g$.
    \item As mentioned above, the term $f(x_0)^2/2g^2$, which did not contribute to Eq.~(\ref{eq:eom}), becomes influential in determining the next approximate stationary point of the Laplace approximation.
\end{itemize}
The consequences of $f_0/f'_0$ not being small can be estimated readily. Consider the next term in the Taylor expansion,
\begin{equation}
    f(x_0 + \delta) = f_0 + f'_0 \delta + \frac{1}{2} \frac{{f^{''}_0}}{{{f'_0}^{2}}} (f'_0 \delta)^2 .
    \label{eq:taylor2}
\end{equation}
Using this equation, we vary $f(x_0+\delta)^2/2g^2$ with respect to the dimensionless variable $\Delta \equiv f'_0 \delta/g.$ 
Note the power of $g$ in this definition, included so that the condition for validity of the Laplace expansion becomes $\Delta \ll 1.$ The variation gives
\begin{equation}
    \left(\frac{1}{g} f_0 + \Delta + \frac{1}{2} \frac{{f^{''}_0}}{{{f'_0}^{2}}} g\Delta^2\right)\left(1+\frac{{f^{''}_A}}{{{f'_0}^{2}}}g\Delta\right).
\end{equation}
Ignoring the second derivative, this implies $\Delta = - f_0/g.$ This is the standard Newton-Raphson update. For the validity of the Laplace approximation, $\Delta \ll 1,$  so we see again that $f_0/g$ is the effective expansion parameter. However, if we include the second derivative terms and estimate $\Delta$ from the resulting quadratic equation, we find that 
\begin{equation}
   \Delta \approx -  \frac{1}{g} f_0 {-} \frac{1}{2g} f_0^2 \frac{{f^{''}_0}}{{{f'_0}^{2}}},
   \label{eq:secondorder}
\end{equation}
which implies that including nonlinearities, with higher derivative terms in the function $f$ for example, cannot be accommodated in the Laplace approximation for arbitrarily small $f_0/g$ in a consistent manner.

\subsection{Correcting the Newton-Raphson method}
Given the inconsistency in the Newton-Raphson method, an improved effective version of Eq.~(\ref{eq:secondorder}) is needed that gives a consistent order by order expansion in $f_0/g,$ but still gives exactly the same formal stationary point equation for $\Delta.$ 
We begin with a clue by observing
\begin{equation}
    \frac{1}{g^2} f(x_0 - f_0/f'_0) \approx \frac{1}{2} \frac{f_0^2}{g^2} \frac{{f^{''}_0}}{{{f'_0}^{2}}} \propto  \bigg(\frac{f_0}{g} \bigg)^2.
\end{equation}
This motivates us to formulate an effective stationary point equation, 
\begin{align}
   \Delta &= -  \frac{1}{g} f_0 - \frac{1}{g^2} f(x_0 + g\Delta/f'_0) \nonumber \\
   &\approx -  \frac{1}{g} f_0 - \frac{1}{g^2} f(x_0 - f_0/f'_0) + \ldots.
   \label{eq:xieom}
\end{align}
{The presence of the second nonlinear term with the lower power of $g$ now becomes influential in determining the consistency of the Laplace approximation, a role absent in the contribution of the term $f(x_0)^2/2g^2$ in Eq.~(\ref{eq:eom}).
Another insight into the lower power of $g$ is that as $g$ increases, $f_0/g$ diminishes, aligning with the strong coupling limit in physics parlance~\citep{svaiter2005strong}, consequently reducing the impact of the additional term in Eq.~(\ref{eq:xieom}).
This term which becomes more important as $g\downarrow 0$ is enforcing the initial ultralocal exact equation: $f(x) = 0.$ The original Newton-Raphson linear term is more important in the strong coupling limit and because it incorporates the derivative of $f,$ and not just the value of $f,$ it corresponds to a `hopping approximation' to the actual equation of motion.
In the Newton-Raphson scenario, the Laplace approximation necessitates further corrections beyond the linear expansion of $f(x)^2/2g^2$ in integral evaluation.
Hence, the inclusion of the additional term with the lower power aims to balance the update of the stationary point between the weak and strong coupling evaluations of the integral.}

To see that the power of $g$ is now exactly correct in a formal expansion, we use the linear Taylor expansion and get 
\begin{equation}
    \frac{1}{g^2} f(x_0 + g\Delta/f'_0) \approx \frac{1}{g^2} (f_0 + g\Delta)
\end{equation}
and inserting this expression back into Eq.~(\ref{eq:xieom}), we get
\begin{equation}
    \Delta = - \frac{1}{g}\left(1+\frac{1}{g}\right) f_0 - \frac{1}{g}\Delta \quad \Leftrightarrow \quad \Delta = -  \frac{1}{g} f_0.
\end{equation}
In other words, the effective equation of motion, when $g$ is small, recapitulates Eq.~(\ref{eq:eom}). 
Finally, we accomplish both objectives of (i) having a leading correction singular in a power of $g$ and (ii) ensuring the term $f(x_0)^2/2g^2$ plays a significant role in the Laplace approximation.
By substituting $\Delta = f'_0 \delta/g$ with $\delta = x - x_0$, the effective stationary point equation of Eq.~(\ref{eq:xieom}) can be expressed as follows:
\begin{eqnarray}
\label{eq:seq0}
\frac{1}{g} f(x) + f_0 + f'_0 (x - x_0) = 0.
\end{eqnarray}

This stationary equation can also be understood more intuitively.
Consider the modified form of the partition function in Eq.~(\ref{eq:Z}) as follows:
\begin{equation}
Z = \int \dd x \exp\left[-\frac{1}{2g^2} \bigg(\alpha f(x) + (1-\alpha) f_L(x) \bigg)^2\right],
\label{eq:newZ}
\end{equation}
where $f_L(x) \equiv f_0 + f'_0(x-x_0)$ represents the linear Taylor approximation of $f(x)$ at $x=x_0$.
This interpolation between $f(x)$ and $f_L(x)$ encapsulates three essential properties:
\begin{itemize}
    \item Setting $\alpha = 1$ recovers the original partition function.
    \item It closely resembles the original partition function when $f_L(x) \approx f(x)$ is a valid approximation.
    \item The term $\alpha f(x)$ serves as a ``tadpole'' that diminishes as $x$ approaches $x_\root$.
\end{itemize}
With a few algebraic manipulations and redefined parameters in Eq.~(\ref{eq:newZ}), we obtain
\begin{equation}
Z = \int \dd x \exp\left[-\frac{1}{2\tilde{g}^2} \bigg(\beta f(x) + f_L(x) \bigg)^2\right],
\label{eq:newZR}
\end{equation}
where $\tilde{g} \equiv g(1-\alpha)$ and $\beta= \alpha/(1-\alpha)$.
Subsequently, the stationary equation leads to:
\begin{eqnarray}
\label{eq:seq}
\beta f(x) + f_0 + f'_0 (x - x_0) = 0.
\end{eqnarray}
Aside from coefficient differences, this equation mirrors Eq.~(\ref{eq:seq0}).
As we will demonstrate, the coefficient acts as a control parameter, rendering this distinction insignificant.

\section{Iterative root-finding}
\label{sec:iterative}

The revised stationary Eq.~(\ref{eq:seq}) serves as an iterative root-finding method. This section elucidates its enhanced convergence, robustness, and interpretation as an iterative algorithm.

\subsection{Iterative algorithm}
We rewrite the stationary Eq.~(\ref{eq:seq}):
\begin{equation}
    \beta f(x(\beta)) + f_0 + f'_0 (x(\beta) - x_0) = 0
    \label{eq:beta}
\end{equation}
to find $x(\beta),$ in place of Eq.~(\ref{eq:eom}). Once we get to Eq.~(\ref{eq:beta}), we can completely ignore the origins of this equation so the numerical coefficient that we have not determined can simply be absorbed into the definition of the abstract parameter $\beta.$

\begin{figure}[t!]
    \centering
    \includegraphics[width=1\columnwidth]{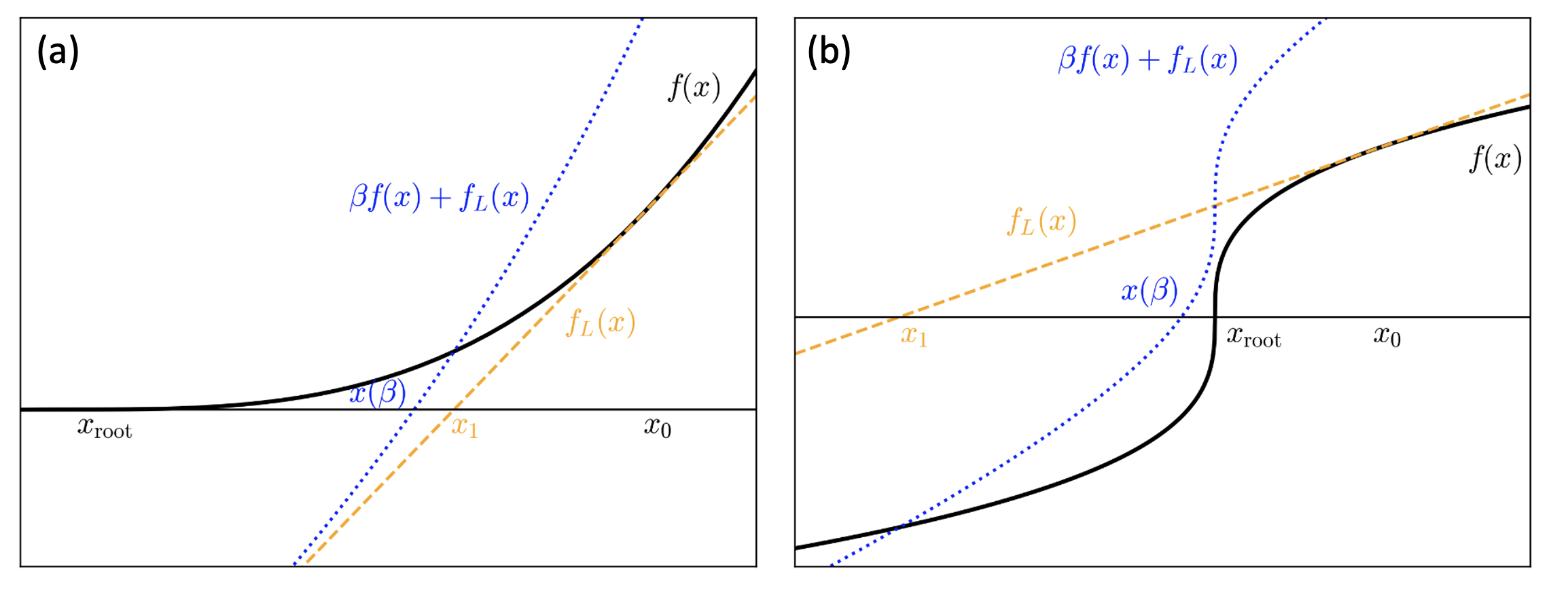}
    \caption{Improved convergence to roots: Comparison between Newton-Raphson method and the revised method. The Newton-Raphson method progresses from an initial point $x_0$ to $x_1$, determined as the solution of $f_L(x_1) = f_0 + f'_0 (x_1 - x_0)=0$, whereas the revised method advances $x_0$ to $x(\beta)$, found as the solution of $\beta f(x) + f_L(x).$ The separation $|x(\beta) - x_0|$ adjusts, either increasing or decreasing, to enhance convergence towards $x_\root$, contingent upon the relationship between $f(x(\beta))$ and $f_0$: (a) when they share the same sign, or (b) when they have opposite signs.}
    \label{fig:fig1}
\end{figure}

Here, we consider intuitively how this new approach works.
Equation~(\ref{eq:beta}) incorporates the complete function $f,$ not just the linear Taylor expansion approximation and therefore makes sense even when $\delta = x(\beta)-x_0$ is not small since, for $\beta$ large, it dominates the other terms in Eq.~(\ref{eq:beta}) to move $\delta$ towards $x_\root,$ as we show geometrically below.
If $f(x(\beta))$ and $f_0$ have the same sign, increasing $\beta$ will result in $|x(\beta) - x_0|$ becoming larger, while the opposite sign will lead to a smaller separation between $x_0$ and $x(\beta).$ This will help the iterative process ($x_0 \rightarrow x(\beta)$) converge faster (Fig.~\ref{fig:fig1}).
For any value $\beta \not= - 1,$ if there is a fixed point, $x^*,$ with a non-singular gradient, we have
\begin{equation}
(\beta + 1) f(x^*) = 0 \ \Leftrightarrow \ f(x^*) = 0.
\end{equation}
The intermediate values of $x(\beta)$ may vary with $\beta$, potentially causing shifts in the basins of attraction as $\beta$ changes (refer to Section~\ref{entropy} for further elucidation). Nonetheless, the fixed point remains a root of $f$, unaffected by $\beta$, a fact that can also be demonstrated through cohomological quantum field theory (Appendix \ref{topological}).

We now explicitly present the revised iterative method as an alternative to the Newton-Raphson method.
Using Eq.~(\ref{eq:beta}), let $x(\beta)$ be defined as the solution of
\begin{equation}
 x(\beta) = N(x_0) - \frac{\beta f(x(\beta))}{f'_0}
\label{eq:basic}
 \end{equation}
with a nonlinear mapping of $N(x) \equiv x - f(x)/f'(x).$ 
At a fixed point, $f(x^*) = 0,$ all $\beta$ dependence vanishes. As we expect based on our discussion above, the fixed point is a root of $f$: $x^* = x_\root.$ 
At $\beta = 0,$ the infinite temperature limit, this is clearly the Newton-Raphson iteration, and does not have any $x(\beta)$ on the right hand side. For $\beta > 0,$ it could also be solved by inserting this definition of $x(\beta)$ into the right hand side repeatedly. At the first order, we get
\begin{equation}
x(\beta) = N(x_0) - \frac{\beta f(N(x_0))}{f'_0}.
\label{eq:iteration0}
\end{equation}
A fixed point of this equation does not immediately imply that the fixed point is a root of $f$, because as it stands this is not the original Eq.~(\ref{eq:basic}). A fixed point here implies only that $f(x^*) + \beta f(N(x^*)) = 0.$ Now this equation is satisfied when $x^*$ is a root, $x_\root$, of $f,$ but the reverse implication is not necessarily true, unless this continues to hold as we vary $\beta$ continuously.

\begin{figure}[t!]
    \centering
    \includegraphics[width=1\columnwidth]{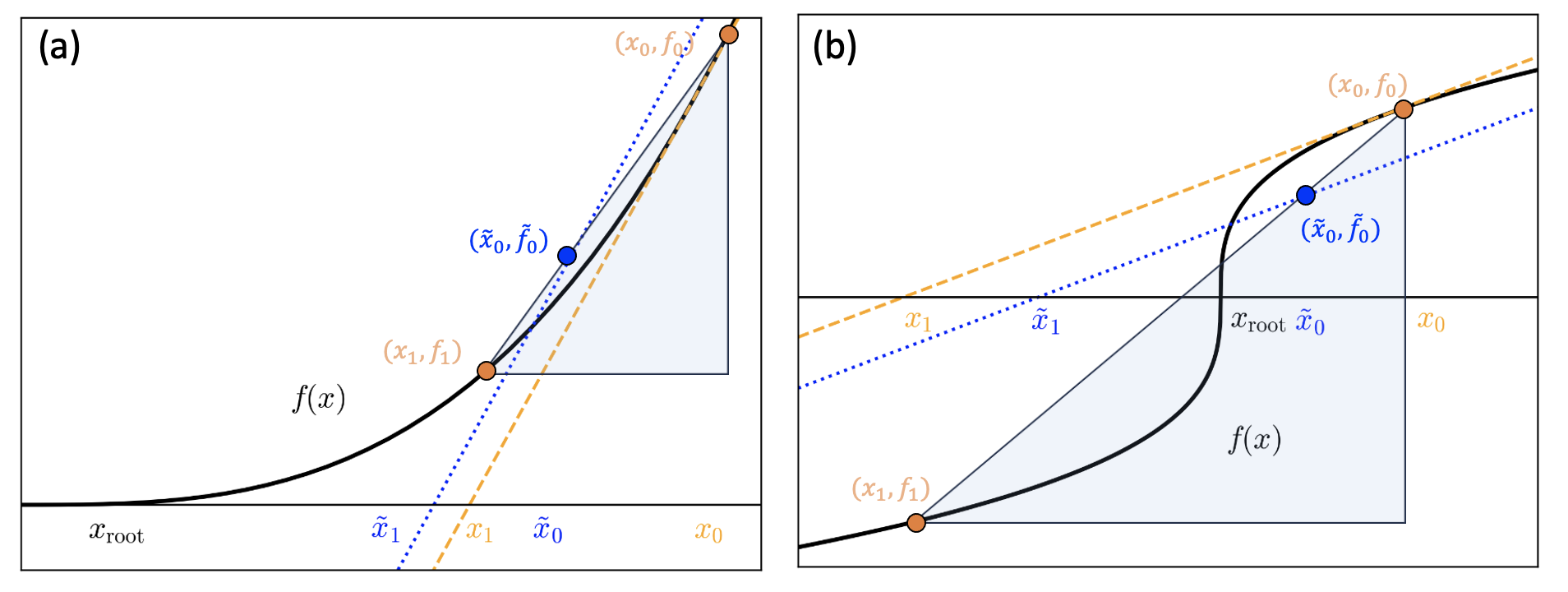}
    \caption{Improved iteration to roots: Comparison between Newton-Raphson method and the revised method. In the Newton-Raphson method, the progression from an initial point $x_0$ to $x_1 = x_0 - f_0/f'_0$ with a slope $f'_0.$ The revised method can be interpreted to provide an interpolated point, $(\tilde{x}_0, \tilde{f}_0)$, between $(x_0, f_0)$ and $(x_1, f_1)$, where $f_i \equiv f(x_i).$ Then, the next iterated point $\tilde{x}_1$ is determined by a linear function crossing $(\tilde{x}_0, \tilde{f}_0)$ with the same slope $f'_0.$ The convergence behavior depends on the relationship between $f_0$ and $f_1$: (a) when they share the same sign, or (b) when they have opposite signs.}
    \label{fig:fig2}
\end{figure}

Let us rewrite the first order equation as
\begin{align}
   x(\beta) &= x_1 - \beta {f_1}/{f'_0} \nonumber \\
   &= (1-\beta)\big(x_0 - {f_0}/{f'_0} \big) + \beta x_1 - \beta {f_1}/{f'_0} \nonumber \\
   &= \tilde{x}_0 - {\tilde{f}_0}/{f'_0},
\end{align}
where we denote $x_1 \equiv N(x_0)= x_0 - f_0/f'_0$ and $f_1 \equiv f(x_1)$, and define $\tilde{x}_0 \equiv (1-\beta)x_0 + \beta x_1$ and $\tilde{f}_0 \equiv (1-\beta)f_0 + \beta f_1$.
Comparing this update with the Newton-Raphson form in Fig.~\ref{fig:fig2}, we see that the function pair $(x_0, f_0)$ determining the next estimate is replaced in our approach with the pair $(\tilde{x}_0, \tilde{f}_0),$ but with the same slope, $f'_0.$ In other words, the line going from $(x_0, f_0)$  to $(x_1, 0)$ with slope $f'_0$ is translated to a parallel line going from $(\tilde{x}_0, \tilde{f}_0),$ to $(x(\beta),0)$ to find the next estimate of the root, $x(\beta).$ Geometrically, this makes the next step in the iteration larger or smaller depending on the relative signs of $f_0$ and $f(x(\beta)),$ as discussed above. 

Were we to repeat this process, we would get
\begin{equation}
x(\beta) = N(x_0) - \frac{\beta f\Big[N(x_0)- \frac{\beta f(N(x_0))}{f'_0} \Big]}{f'_0}.
\label{eq:iteration}
\end{equation}
and so on. Each order in $\beta$ has one more function evaluation and is higher order in $\beta.$ It turns out that either of these definitions of iterative determination of $x(\beta)$ gives exactly the same fixed point structure, which exhibits a fractal nature.
A fixed point value that is a fixed point of both these versions will in fact be a root.
Henceforth, we proceed with the simplest choice, which can be described by two iterative steps:
\begin{align}
\hat{x}_{n+1} &= N(x_n) = x_n  - \frac{f({x}_{n})}{f'(x_n)}, \nonumber \\
x_{n+1} &= \hat{x}_{n+1} - \beta \frac{f(\hat{x}_{n+1})}{f'(x_n)}.
\label{eq:iteration_beta}
\end{align}
This completes our derivation of Eq.~(\ref{eq:newnewton}).
It is noteworthy that this modified Newton-Raphson method, with a specific value of $\beta=1$, has been derived in \citep{amat2004dynamics, amat2010chaotic, Chun2005}.
The derivation in \citep{Chun2005} based on the Adomian method will be further discussed in the subsequent section.

\subsection{Convergence} \label{faster}
We now demonstrate that there is actually a benefit of this $\beta$ dependent reformulation in terms of root-finding performance. 
The quadratic convergence of the Newton-Raphson method, near a root, $x_\root,$ is demonstrated by Taylor expanding 
\begin{equation}
f(x) \approx  f'_\root(x-x_\root) +\frac{1}{2}f^{''}_\root (x-x_\root)^2,
\end{equation}
where $f'_\root \equiv f'(x_\root)$ and $f^{''}_\root \equiv f^{''}(x_\root)$.
By inserting this Taylor approximation into the Newton-Raphson update of $x_{n+1} = N(x_n)$, 
we obtain the the expected quadratic convergence: 
\begin{equation}
 x_{n+1} -x_\root \approx \frac{1}{2} \frac{{{f^{''}_\root }}}
{f'_\root +f^{''}_\root (x_n-x_\root)} (x_n-x_\root)^2.
\end{equation}
Similarly, for our $\beta$ dependent update, we get
\begin{equation} \label{eq:update_order}
 x_{n+1} -x_\root \approx \frac{1-\beta}{2} \frac{f^{''}_\root}{f'_\root}(x_n-x_\root)^2,
\end{equation}
exactly in line with the intuition we gave above. Of course, the number of iterations is never zero so, more precisely, this only shows that the convergence at $\beta = 1$ is faster than quadratic. See \citet{Homeier2004, Homeier2005}\ for the significance of cubic convergence in modified Newton-Raphson methods. 
Here we see that for problems where the first derivative at the root is finite, we should use the largest value of $\beta$ for which a fixed point exists. However, as we shall show explicitly in the case of $f(x) = x^{1/3},$ a fixed point may not exist for $\beta = 1$.

Close to a root, the order of convergence depends on the parameter $\beta$, being quadratic for the Newton-Raphson case with $\beta = 0$. Theorem 1 in \cite{grau2004} proves that the order of convergence is maximized for $\beta = 1$, achieving local order three for a similar iterative algorithm.
The cubic convergence for $\beta=1$ has also been demonstrated in other studies~\citep{traub1982iterative, petkovic2014multipoint}.
However, the theorem does not specify the convergence order for other values of $\beta$.
We conducted numerical validation to assess the convergence order across a range of nonlinear functions, and confirmed the quadratic convergence for $\beta \neq 1$ and the cubic convergence only for $\beta = 1$ (refer to Section~\ref{sec:results} in details).

\begin{table*}[!t]
\begin{tabular}{p{5.5cm}|>{\centering\arraybackslash}p{0.9cm}>{\centering\arraybackslash}p{0.65cm}>{\centering\arraybackslash}p{0.65cm}>{\centering\arraybackslash}p{0.65cm}>{\centering\arraybackslash}p{0.65cm}|>{\centering\arraybackslash}p{0.9cm}>{\centering\arraybackslash}p{0.6cm}>{\centering\arraybackslash}p{0.6cm}>{\centering\arraybackslash}p{0.6cm}>{\centering\arraybackslash}p{0.6cm}|>{\centering\arraybackslash}p{0.9cm}>{\centering\arraybackslash}p{0.6cm}>{\centering\arraybackslash}p{0.5cm}>{\centering\arraybackslash}p{0.6cm}>{\centering\arraybackslash}p{0.6cm}}
    \hline
{\multirow{2}{*}{Functions}}  &   \multicolumn{5}{c|}{Iteration number} &   \multicolumn{5}{c|}{Convergence percentage} &   \multicolumn{5}{c}{Computation time}  \\
    \cline{2-16}
& {\footnotesize$\beta=$}$-1$ & $-0.5$ & $0$ & $0.5$ & $1$ & {\footnotesize$\beta=$}$-1$ & $-0.5$ & $0$ & $0.5$ & $1$ & {\footnotesize$\beta=$}$-1$ & $-0.5$ & $0$ & $0.5$ & $1$ \\
    \hline
{\footnotesize $f_1(x) = (x^2 - 1)(x^2 + 1)$}  & 18.4  & 15.7 & 11.3 & 11.0 & \bf{10.5} & 75  & 80 & \bf{100} & \hfil 91 & \hfil 90 & 1.50 & 1.29 & 1.0 & 0.92 & \bf{0.89}   \\
    \hline
{\footnotesize $f_2(x) = x^3 - 1$} & 15.5 & 13.3 & \hfil 9.1 & \hfil 9.1 & \hfil \bf{8.4} & 92 & 96 & \bf{100} & \hfil 99 & \hfil 99 & 1.55 & 1.35 & 1.0 & 0.94 & \bf{0.93} \\
\hline
{\footnotesize $f_3(x) = x^{12} - 1$} & 18.5 & 16.3 & 16.4 & 12.7 & \bf{11.5} & 44 & 56 & \hfil \bf{84} & \hfil 67 & \hfil 66 & 1.11 & 0.93 & 1.0 & 0.73 & \bf{0.67} \\
\hline
{\footnotesize $f_4(x) = (x^2 - 4)(x + 1.5)(x - 0.5)$} & 12.6 & 10.4 & \hfil 8.0 & \hfil 7.4 & \hfil \bf{6.4} & 98 & 97 & \bf{100} & \bf{100} & \bf{100} & 1.50 & 1.27 & 1.0 & 0.89 & \bf{0.79} \\
\hline
{\footnotesize $f_5(x) =(x+2)(x+1.5)^2 (x-0.5)(x-2)$} & 12.1 & \bf{10.3} & 28.9 & 25.6 & 22.2 & 42 & 43 & \hfil 96 & \hfil \bf{99} & \hfil \bf{99} & 0.42 & \bf{0.35} & 1.0 & 0.85 & 0.76 \\
\hline
{\footnotesize $f_6(x) = \sin(x)$} & \hfil 9.7 & \hfil 7.7 & \hfil 6.6 & \hfil 6.0 & \hfil \bf{5.5} & 87 & 84 & \bf{100} & \hfil 95 & 
 \hfil 92 & 1.42 & 1.09 & 1.0 & 0.87 & \bf{0.82} \\
\hline
{\footnotesize $f_7(x) = (x - 1)^3 + 4(x-1)^2 - 10$} & 17.9 & 14.5 & \hfil 9.7 & \hfil 9.1 & \hfil \bf{8.4} & 97 & 97 & \bf{100} & \bf{100} & \bf{100} & 1.74 & 1.44 & 1.0 & 0.92 & \bf{0.83} \\
\hline
{\footnotesize $f_8(x) = \sin(x-14/10)^2 - (x - 14/10)^2 + 1$} & 11.3 & \hfil 9.4 & \hfil 8.9 & \hfil 7.5 & \hfil \bf{6.5} & 60 & 59 & \hfil \bf{71} & \hfil 64 & \hfil 64 & 1.25 & 1.04 & 1.0 & 0.86 & \bf{0.76} \\
\hline
{\footnotesize $f_9(x) = x^2 - e^x - 3x + 2$} & \hfil 8.1 & \hfil 7.1 & \hfil 6.2 & \hfil 5.8 & \hfil \bf{5.0} & 97 & 97 & \bf{100} & \hfil 99 & \hfil 98 & 1.20 & 1.13 & 1.0 & 0.93 & \bf{0.79} \\
\hline
{\footnotesize $f_{10}(x) = \cos(x-3/4) - x + 3/4$} & 10.8 & \hfil 8.8 & \hfil 8.6 & \hfil 7.1 & \hfil \bf{6.3} & 61 & 63 & \bf{93} & \hfil 75 & \hfil 72 & 1.23 & 0.97 & 1.0 & 0.8 & \bf{0.72} \\
\hline 
{\footnotesize $f_{11}(x) = (x + 1)^3 - 1$} & 15.1 & 12.7 & \hfil 9.0 & \hfil 8.9 & \hfil \bf{8.1} & 93 & 96 & \bf{100} & \hfil 99 & \hfil 99 & 1.50 & 1.29 & 1.0 & 0.96 & \bf{0.85} \\
\hline
{\footnotesize $f_{12}(x) = (x-2)^3 - 10$} & 17.5 & 14.7 & \hfil 9.3 & \hfil 9.8 & \hfil \bf{8.9} & 87 & 92 & \bf{100} & \hfil 97 & \hfil 98 & 1.84 & 1.54 & 1.0 & 1.01 & \bf{0.96} \\
\hline
{\footnotesize $f_{13}(x) = (x+5/4) e^{(x+5/4)^2} - \sin(x+5/4)^2 + 3\cos(x+5/4) + 5$} & 17.0 & 13.2 & 11.6 & \hfil 9.7 & \hfil \bf{9.4} & 55 & 60 & \bf{96} & \hfil 82 & \hfil 80 & 1.46 & 1.12 & 1.0 & 0.81 & \bf{0.80} \\
\hline
{\footnotesize $f_{14}(x) = x + \sin(2/x)  x^2$} & 15.5 & 13.5 & 10.9 & \hfil \bf{9.5} & 10.3 & 35 & 57 & \bf{99} & \hfil 88 & \hfil 74 & 1.41 & 1.22 & 1.0 & \bf{0.86} & 1.08 \\
\hline
\end{tabular}
\caption{\label{table:1}
Numerical performance of the extended Newton-Raphson method. For five different values of $\beta$, we evaluate: (i) the average number of iterations required for convergence (Iteration number); (ii) the percentage of initial points that converge (Convergence percentage); and (iii) the relative computation time compared to the $\beta=0$ case (Computation time). Note that $\beta=0$ represents the original Newton-Raphson method. In each category, the best performance—shortest iteration number, highest convergence percentage, and least computation time among different $\beta$ values—is highlighted in bold.
}
\end{table*}

\subsection{Cube root function} \label{cuberoot}
The cube root is a well-known analytically solvable function that is not well-suited for finding roots using the Newton-Raphson algorithm.
In fact, for $f(x) = x^{1/3},$ we get $f(x)/f'(x) = 3x.$ It follows that 
\begin{equation}
x_{n+1} = N(x_n) = -2 x_{n}, 
\end{equation}
which clearly does not converge. For our $\beta$ dependent update, Eq.~(\ref{eq:iteration_beta}), we have
$N(x) = -2x$, and so $f(N(x))/f'(x) = (-2x)^{1/3}/(x^{-2/3}/3) = -2^{1/3}\cdot 3x.$ Therefore, Equation~(\ref{eq:iteration_beta}) implies
\begin{equation}
x_{n+1} = x_n \left(-2 + 3\cdot 2^{1/3}\beta\right).
\end{equation}
For convergence we must have 
\begin{equation}
|-2 + 3\cdot 2^{1/3}\beta| < 1,
\end{equation}
which requires 
\begin{equation}
\frac{1}{3 \cdot 2^{1/3}} < \beta < \frac{1}{2^{1/3}}.
\end{equation}
Moreover, the swiftest convergence occurs at $\beta_\mathrm{min} = 2^{2/3}/3.$ It is worth noting that $\beta = 0$ (Newton-Raphson method) and $\beta = 1$ (Adomian's method) do not fall within the range of convergent values. These observations regarding the strong-coupling approach to root determination of the cube root function are empirically validated in one of the examples. Remarkably, by tuning $\beta$ to this particular value, convergence becomes seemingly independent of the initial point in this instance.
This surprising convergence is indicative of a more interesting role for $\beta$ in the size of the basins of attraction to a given root as we shall show in Section~\ref{entropy}.

\subsection{\label{sec:results}Numerical results}

We now explicitly compare key numerical figures of the extended Newton-Raphson method. Specifically, the original Newton-Raphson method corresponds to the case of $\beta=0$ in the parameterized Newton-Raphson method. To conduct the comparison, we tested various types of nonlinear functions, some of which were adopted from~\citet{weerakoon2000variant}. 

Tables~\ref{table:1} and \ref{table:2} show four metrics:
\begin{itemize}
\item Iteration number: The average number of iterations per initial point after discarding divergent initial points.
\item Convergence percentage: The percentage of initial points that successfully converge to a root before reaching a preset maximum number of iterations (set to 50 iterations in this study).
\item Computation time: The average execution time per point relative to the Newton-Raphson method. Divergent initial points are excluded from this calculation. 
\item Convergence order: The estimated order of convergence as iteration progresses.
\end{itemize}

Given that finding real roots for real-valued function $f(x)=0$ with $x\in \mathbb{R}$ can be extended to finding zeros for complex-valued functions $f(z)=0$ with $z\in \mathbb{C}$, we explore the Newton-Raphson methods in the complex plane.
For each function, we iterate the method starting from an initial point $z_0 \in \mathbb{C}$, chosen from a regular grid in the square $[-2, 2] \times [-2, 2]$. We select a grid of 1000 by 1000 points evenly spaced across this square and iterate each initial point until the criterion $|z_{n+1} - z_n| < \varepsilon$ is satisfied, with a threshold $\varepsilon = 10^{-14}$. If the iteration does not satisfy the criterion within 50 iterations, we consider the trajectory divergent.
Here, the convergence order was estimated using the following formula \cite{grau2010}:
\begin{equation}\label{eq:ACOC}
q_n = \frac{\log|e_{n+1}/e_n|}{\log|e_n/e_{n-1}|},
\end{equation}
where $e_n \equiv z_n - z_{n-1}$.
When $|e_n| < \varepsilon$ the iterations stop and we pick the last value of the series $q_n$ as the order of convergence.
For each function, we select an initial point that converges to one of the roots for at least eight iterations to obtain a valid estimation of the convergence.

We investigate the impact of the parameter $\beta$ in the extended Newton-Raphson method by using some representative values of $\beta \in \{-1, -0.5, 0 , 0.5, 1\}$.
We found that positive values of $\beta$ improve the efficiency of the root-finding iteration compared to the original Newton-Raphson method ($\beta = 0$). Specifically, the case $\beta = 1$ consistently shows superior performance (Tab.~\ref{table:1}).
The mean number of iterations decreases monotonically as $\beta$ changes from $-1$ to $1$. However, this trend for positive $\beta$ is accompanied by higher instabilities as the convergence percentage of initial points decreases. The relative computation time compared to the original Newton-Raphson method is reduced. Since additional operations for the extended method are negligible, the reduced number of iterations directly leads to shorter computation times.

In general, negative values of $\beta$ lead to poorer performance, except for the function $f_5$. A closer examination of the convergence percentage reveals more than half of the initial points do not converge for negative $\beta$. This can be understood visually from Fig.~\ref{fig:fig1}, where unlike postive $\beta$, negative $\beta$ causes the next iteration points to move farther away from the roots.

Table~\ref{table:2} showed quadratic convergence for $\beta = 0$ and cubic convergence for $\beta = 1$. Notably, we encountered anomalous convergence in the case of the function $f_6(x) = \sin(x)$, which exhibits vanishing curvature with $f_6^{''}(x_\root)=0$. Additionally, we noted linear convergence near a doubly degenerate real root for $f_5(x)$.

\begin{figure*}[t!]
    \centering
    \includegraphics[width=1\columnwidth]{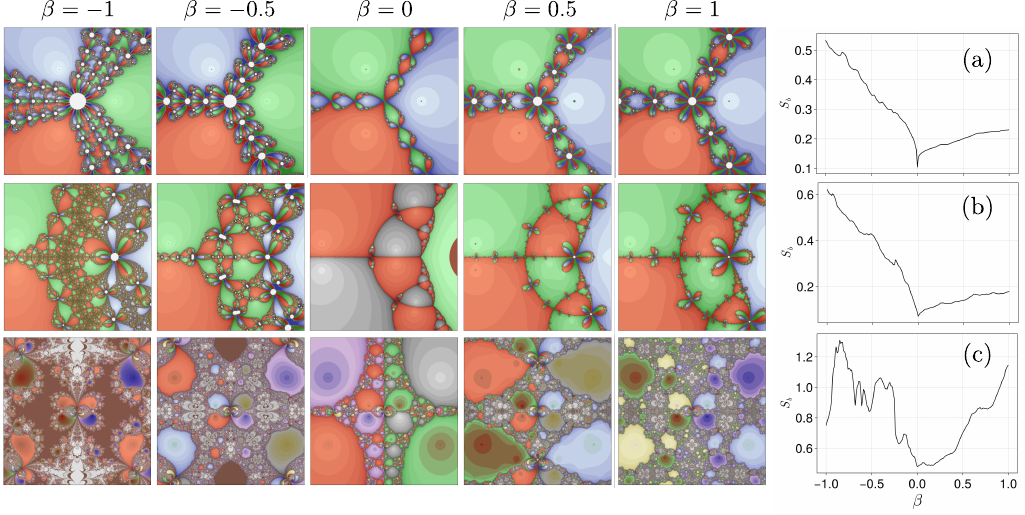}
    \caption{Newton fractals and basin entropy. Colors represent different basins of attractors (roots). In particular, the colors in the basins of attraction have been shaded such that darker colors correspond to longer iterations to converge to roots. The basins of three functions are computed for $\beta \in \{ -1 , -0.5, 0 , 0.5 , 1\}$: (a) top row $f_2(z) = z^3 -1$, (b) central row $f_7(z) = (z - 1)^3 + 4(z-1)^2 - 10$, and (c) bottom row $f_{14}(z) = z + z^2 \sin(2/z)$. The corresponding basin entropy $S_b$ in the right plot quantifies the unpredictability of roots as a function of $\beta$ over the range $[-1,1]$ in steps of $0.01$. $S_b$ has been computed on a grid of $1000$ times $1000$ initial points with a covering of boxes of size $20$ times $20$. The method to determine the final root is described in Sec.~\ref{sec:results}.}
    \label{fig:frac_bas2}
\end{figure*}

\subsection{Basin entropy}
\label{entropy}

For nonlinear functions with multiple roots, different initial points $z_0$ approach different roots or zeros or attractors. Consequently, the set of initial points converging to the same root $z^*$ is defined as the basin of attraction of the root $z^*$. It is well known that the basin of attractors exhibits a fractal structure for the Newton-Raphson method. Here, we investigate whether our iteration method alters the Newton-Raphson basin of attractors.
A remarkable characteristic of the basins resulting from the numerical method is the intricate nature of the boundary between them. For initial points located on this boundary, the uncertainty in the final root is maximized.
This aspect could be perceived as a drawback of the method, as a ``poor'' initial point might result in convergence to any root of the function. Nonetheless, for functions with unknown roots, our interest lies in identifying all roots in principle. As we will explore in the next section, this characteristic can also contribute to an enhancement in the speed of convergence.

We investigate the impact of our method on the boundaries of the basins generated by iterations from initial points. Upon initial inspection of the basins for the function $f_2(z) = z^3 - 1$, a significant alteration in the boundary between $\beta = -1$ and $\beta = 1$ becomes evident (Fig.~\ref{fig:frac_bas2}(a) and top row). Formally, in the complex plane, this boundary constitutes the set of points that never converge to a root and forms a Julia set~\citep{falconer2004}. Furthermore, the example illustrated in Fig.~\ref{fig:frac_bas2}(a) top row exhibits the property of Wada, wherein a single boundary simultaneously separates three basins. This property entails a unique form of unpredictability, as a point on the boundary can ultimately converge to any of the roots. Similar phenomena have been observed on boundaries for other modified Newton methods dependent on parameters, as discussed in \cite{susanto2009}, where the modifications introduce non-trivial transformations of the boundaries. Here, we provide quantitative insight into these transformations.

The unpredictability of the final root for initial points near the boundary can be quantified using the \textit{basin entropy}~\citep{Daza2016}. This metric assesses the local entropy within boxes of side length $l$ by initially estimating the probability $p_i$ of each final attractor $i$ inside the box with the na\"\i ve frequency estimator. These estimated probabilities serves to compute the Gibbs (or Shannon) entropy of the box: $S_b = -\sum_i p_i\log p_i$. The basin entropy is then calculated as the average of the box entropy with a covering of non-overlapping boxes over the portion of the phase space studied.

In the examples shown in Fig.~\ref{fig:frac_bas2}, the  basins of the roots are computed over a grid of $1000\times 1000$ initial points in the complex plane using the methods of Sec.~\ref{sec:results}. A box size of $20 \times 20$ initial points is used to estimate the $p_i$ probabilities. 

In Fig.~\ref{fig:frac_bas2}, we depict the evolution of the basin entropy for three functions across the range of $\beta$ values from -1 to 1. It is evident that unpredictability consistently increases for $\beta\neq 0$. While we cannot offer a rigorous explanation for this phenomenon, we can provide heuristic arguments regarding the emergence of new fractal structures in the phase space. 
One mechanism contributing to the formation of fractal boundaries is the stretching and folding action of small areas. The roots of the derivative $f'(x)$ serve as sources of instabilities in the Newton-Raphson method, as initial points near these points can be dramatically dispersed. Subsequent iterations tend to bring these points back towards one of the roots. Small variations in the initial points near these roots can result in convergence to any of the roots, characteristic of chaotic behavior.

We consider our new iterative method, expressed as $x_{n+1} = N(x_n) - \beta {f(N(x_n))}/{f'(x_n)}$, where $N(x) = x - f(x)/f'(x)$. When close to a singular $x_\singular$, where $f'(x_\singular)=0$, we can approximate a linear relation, $f'(x) \simeq c_2 \delta x$ with $c_2 \equiv f^{''}(x_\singular)$ and $\delta x \equiv x - x_\singular$. Consequently, $N(x) \simeq - c_0/\delta x$ with $c_0 \equiv f(x_\singular)$. The amplitude of the next iteration, given by
\begin{equation}
    x_{n+1} \simeq - \frac{c_0}{\delta x_n} - \frac{\beta f(-c_0/\delta x_n)}{c_2 \delta x_n}, 
\end{equation} 
depends on the proximity to the singularity, $\delta x_n = x_n - x_\singular$, and the value of $f(-c_0/\delta x_n)$. Generally, uncertainty increases near these points due to these significant jumps. For complex functions, new structures emerge near the singularity, visible as the blobs in Figs.~\ref{fig:frac_bas2} for $\beta \neq 0$.

This rise in basin entropy signifies an enhanced mixing property in the phase space, facilitating more thorough exploration of roots. The algorithm can traverse a broader range of regions in the phase space before settling on a local solution. A recent study implemented a deflated version of the Newton-Raphson method~\citep{Cisternas2024}, aiming to reveal additional roots of a function by avoiding convergence to already known solutions. Similarly, our parameterized Newton-Raphson method enhances root searching efficiency in phase space due to the fractal nature of the boundary.

\section{Link to the Adomian method}
\label{sec:adomian}
The self-consistent stationary Eq.~(\ref{eq:basic}) is reminiscent of the Adomian method \citep{adomian2013solving, rach2012bibliography}. Notably, the Adomian method has previously yielded an enhanced Newton-Raphson approach, corresponding precisely to the scenario when $\beta=1$ in our formulation. In this section, we provide an overview of the Adomian method, establish connections between our approach and the Adomian method, and extend the fixed $\beta$ method to an annealing approach with varying $\beta$.

\subsection{Adomian method}
 We first give an extensive introduction to the Adomian method because it may not be familiar to all readers.
The canonical form for employing the Adomian decomposition is given by
\begin{equation}
a = C + F(a),
\label{eq:adobasic}
\end{equation}
where $C$ is a constant and $F(a)$ is a nonlinear function of the variable $a$.
Utilizing a series solution approach, we express:
\begin{align}
a &= a_0 + a_1 + a_2 + \cdots, \nonumber \\
F(a) &= A_0(a_0) + A_1(a_0, a_1) + A_2(a_0, a_1, a_2) + \cdots.
\label{eq:series}
\end{align}
Equating terms between the two sides in Eq.~(\ref{eq:adobasic}) requires explicit choices because there is no formal expansion parameter in the series. With specific choices, this leads to the following relationships:
\begin{align}
a_0 &= C, \nonumber \\
a_1 &= A_0(a_0), \nonumber \\
a_2 &= A_1(a_0, a_1), \nonumber \\
&\cdots
\label{eq:arbitrary}
\end{align}
where $A_n$ represents the Adomian polynomial:
\begin{equation}
A_n = \frac{1}{n!} \left[\frac{d^n}{d \beta^n} F(a_0 + \beta a_1 + \beta^2 a_2 + \cdots) \right]_{\beta = 0}.
\label{eq:adomian_polynomial}
\end{equation}
Here are a few initial terms:
\begin{align}
A_0(a_0) &= F(a_0), \nonumber \\
A_1(a_0, a_1) &= a_1 F'(a_0).
\end{align}
We emphasize that there is no unambiguous way to compare terms between the left and right hand sides of Eq.~(\ref{eq:arbitrary}):
For example, 
\begin{align}
a_0 &= C - \epsilon, \nonumber \\
a_1 &= \epsilon - \epsilon^2 + A_0(a_0), \nonumber \\
a_2 &= \epsilon^2 - \epsilon^3 + A_1(a_0, a_1), \nonumber \\
&\vdots
\end{align}
for an arbitrary $\epsilon$ with $|\epsilon| < 1,$ also solves Eq.~(\ref{eq:adobasic}). 
Results on convergence have been proved with hypotheses on the size of derivatives and it has been noted that the decomposition
must be chosen appropriately~\citep{Abbaoui1994}, as there is no canonical choice for the matching of terms.

Let us apply the Adomian decomposition to the root-finding problem, $f(x) = 0$.
Utilizing a linear approximation of $f(x)$ in proximity to $x$ gives:
\begin{equation}
f(x-a) = f(x) - a f'(x).
\end{equation}
Rearranging this equation into the canonical form for Adomian decomposition yields:
\begin{equation}
a = \underbrace{\frac{f(x)}{f'(x)}}_{C} \underbrace{- \frac{f(x-a)}{f'(x)}}_{F(a)}.
\end{equation}
By focusing solely on the zero-th order term, $a_0 = C = f(x)/f'(x)$, we can approximate $a \approx a_0.$ 
The iterative update of $x_{n+1} = x_n - a$ is expected to satisfy $f(x_{n+1}) = f(x_n -a) = 0$ for root-finding. Then, $x_n - x_{n+1} = a = a_0$ results in 
\begin{equation}
x_{n+1} = x_n - \frac{f(x_n)}{f'(x_n)},
\end{equation}
which is equivalent to the Newton-Raphson method.

\citet{Chun2005} has expanded the analysis by considering an additional step:
\begin{equation}
f(x-a) = f(x) - a f'(x) + \underbrace{\frac{1}{2} a^2 f^{''}(x)}_{g(x,a)}.
\label{eq:chun}
\end{equation}
Once again, rearranging this equation, assuming $f(x-a) = 0$ at $x-a$, allows it to be expressed in the canonical form for Adomian decomposition:
\begin{equation}
a = \underbrace{\frac{f(x)}{f'(x)}}_{C} + \underbrace{\frac{g(x,a)}{f'(x)}}_{F(a)}.
\end{equation}
This decomposition is certainly not unique (see \citet{Abbasbandy2003} for an alternative decomposition).

This time, considering $a \approx a_0 + a_1$ up to the first order, the Adomian decomposition provides:
\begin{align}
a_0 &= C = \frac{f(x)}{f'(x)}, \\
a_1 &= F(a_0) = \frac{g(x,a_0)}{f'(x)} = \frac{f(x-a_0)}{f'(x)},
\end{align}
where we utilized $g(x, a_0) = f(x-a_0) - f(x) + a_0 f'(x) = f(x-a_0)$.
Consequently,
\begin{equation}
a \approx a_0 + a_1 = \frac{f(x)}{f'(x)} + \frac{f(x-a_0)}{f'(x)}.
\end{equation}
Again, given that $x_n - x_{n+1} = a = a_0 + a_1$, this leads to 
\begin{align}
x_{n+1} &= x_{n} - \frac{f(x_{n})}{f'(x_{n})} - \frac{f\bigg(x_{n}-\frac{f(x_{n})}{f'(x_{n})}\bigg)}{f'(x_{n})} \nonumber \\
&= N(x_{n}) - \frac{f(N(x_{n}))}{f'(x_{n})},
\end{align}
which corresponds to Eq.~(\ref{eq:iteration}) in our formulation with $\beta=1$.

The Adomian polynomial approach, when truncated to a specific number of terms in $a_i,\ i\le m$, offers a pathway to derive alternative root-finding iterations. Although the higher-order forms involve algebraic complexities, these truncations result in progressively higher orders of convergence for well-behaved functions. This observation led to Chun's conjecture~\citep{Chun2005} that such truncations exhibit convergence order $m+2.$ Therefore, when $m=0$, it corresponds to the Newton-Raphson iteration, which demonstrates quadratic convergence.

Let us revisit our self-consistent stationary Eq.~(\ref{eq:basic}):
\begin{equation}
x(\beta) = \underbrace{N(x_0)}_{C} \underbrace{-\frac{\beta f(x(\beta))}{f'_0}}_{F(x(\beta))}.
\label{eq:hightemp}
\end{equation}
This equation adheres to the canonical form of the Adomian method. One notable distinction is the inclusion of the parameter $\beta.$ 
In the limit of small $\beta$, {the so-called} ``high-temperature'' limit, the second term can be regarded as a perturbation.
We now explore the high-temperature limit, and consider a series solution:
\begin{equation}
    x(\beta) = a_0 + \beta a_1 + \beta^2 a_2 + \ldots.
\end{equation}
Here, the nonlinear function $F(x(\beta))$ can also be expressed using the Adomian polynomial in Eq.~(\ref{eq:adomian_polynomial}). When we truncate $x(\beta) = a_0 + \beta a_1$ to the first order of $\beta$, we obtain:
\begin{align}
    a_0 &= N(x_0) \nonumber \\
    a_1 &= - \frac{f(a_0)}{f'_0}.
\end{align}
This yields:
\begin{equation}
    x(\beta) = a_0 + \beta a_1 = N(x_0) - \frac{\beta f(N(x_0))}{f'_0},
\label{eq:adomian_iteration}
\end{equation}
which corresponds to Eq.~(\ref{eq:iteration0}). Hence, our formulation can be interpreted within the context of the Adomian method. However, the inclusion of the auxiliary parameter $\beta$ holds significant importance in the Adomian decomposition. It ensures an {unambiguous} comparison between series and the Adomian polynomial, providing an exact term-by-term mapping with correct orders of $\beta$.

We note that there is no need to use Eq.~(\ref{eq:series}) at all as we can work directly with $x(\beta)$ in our iteration, Eq.~(\ref{eq:adomian_iteration}). Also noteworthy is the fact that as we directly find $x(\beta)$ fixed points, no derivatives of $f$ beyond the first are ever explicitly required in the iteration.

\begin{table*}[!t]
\setlength\tabcolsep{0.1pt}%
\begin{tabular}{p{5.5cm}|>{\centering\arraybackslash}p{1cm}>{\centering\arraybackslash}p{1cm}>{\centering\arraybackslash}p{1cm}|>{\centering\arraybackslash}p{1cm}>{\centering\arraybackslash}p{1cm}>{\centering\arraybackslash}p{1cm}|>{\centering\arraybackslash}p{1cm}>{\centering\arraybackslash}p{1cm}>{\centering\arraybackslash}p{1cm}|>{\centering\arraybackslash}p{1cm}>{\centering\arraybackslash}p{1cm}>{\centering\arraybackslash}p{1cm}}
    \hline
{\multirow{2}{*}{Functions}}  &   \multicolumn{3}{c|}{Iter. number} &   \multicolumn{3}{c|}{Conv. percentage} &   \multicolumn{3}{c|}{Comp. time} & \multicolumn{3}{c}{Conv. order}  \\
    \cline{2-13}
& {\small$\beta=0$} & {\small$\beta=1$} & {\small$\beta_n$} & {\small$\beta=0$} & {\small$\beta=1$} & {\small$\beta_n$} & {\small$\beta=0$} & {\small$\beta=1$} & {\small$\beta_n$} & {\small$\beta=0$} & {\small$\beta=1$} & {\small$\beta_n$} \\
    \hline
{\footnotesize $f_1(x) = (x^2 - 1)(x^2 + 1)$}  & 11.3  & 10.5 & \hfil \bf{7.9} & \bf{100} & \hfil 90 & \bf{100}  & 1.0 & \bf{0.89} & 1.02 & 2.00 & 2.99 & \bf{4.08}  \\
\hline 
{\footnotesize $f_2(x) = x^3 - 1$}  & \hfil 9.1 & \hfil 8.4 & \hfil \bf{6.5} & \bf{100} & \hfil 99 & \bf{100} & 1.0 & \bf{0.93} & 1.07 & 2.00 & 3.01 & \bf{4.14} \\
\hline
{\footnotesize $f_3(x) = x^{12} - 1$}  & 16.4 & \bf{11.5} & 12.8 & \hfil 84 & \hfil 66 & \hfil \bf{87} & 1.0 & \bf{0.67} & 1.26 & 2.00 & 2.92 & \bf{4.17} \\
\hline
{\footnotesize $f_4(x) = (x^2 - 4)(x + 1.5)(x - 0.5)$}  & \hfil 8.0 & \hfil 6.4 & \hfil \bf{5.4} & \bf{100} & \bf{100} & \bf{100} & 1.0 & \bf{0.79} & 1.07 & 1.99 & 2.94 & \bf{3.95} \\
\hline
{\footnotesize $f_5(x) =(x+2)(x+1.5)^2 (x-0.5)(x-2)$}  & 28.9 & 22.2 & \bf{18.6} & \hfil 96 & \hfil 99 & \bf{100} & 1.0 & \bf{0.76} & 1.17 & 0.99 & \bf{1.00} & \bf{1.00} \\
\hline
{\footnotesize $f_6(x) = \sin(x)$}  & \hfil 6.6 & \hfil 5.5 & \hfil \bf{5.0} & \bf{100} & \hfil 92 & \hfil 99 & 1.0 & \bf{0.82} & 1.18 & 3.00 & 4.97 & \bf{6.85} \\
\hline 
{\footnotesize $f_7(x) = (x - 1)^3 + 4(x-1)^2 - 10$}  & \hfil 9.7 & \hfil 8.4 & \hfil \bf{7.3} & \bf{100} & \bf{100} & \bf{100} & 1.0 & \bf{0.83} & 1.44 & 2.00 & 3.02 & \bf{4.36} \\
\hline
{\footnotesize $f_8(x) = \sin(x-14/10)^2 - (x - 14/10)^2 + 1$}  & \hfil 8.9 & \hfil 6.5 & \hfil \bf{6.1} & \hfil 71 & \hfil 64 & \hfil \bf{73} & 1.0 & \bf{0.76} & 1.27 & 2.00 & \bf{2.98} & 2.58 \\
\hline
{\footnotesize $f_9(x) = x^2 - e^x - 3x + 2$}  & \hfil 6.2 & \hfil 5.0 & \hfil \bf{4.2} & \bf{100} & 98 & \bf{100} & 1.0 & \bf{0.79} & 0.99 & 2.00 & 3.09 & \bf{4.25} \\
\hline
{\footnotesize $f_{10}(x) = \cos(x-3/4) - x + 3/4$}  & \hfil 8.6 & \hfil \bf{6.3} & \hfil \bf{6.3} & \hfil \bf{93} & \hfil 72 & \hfil \bf{93} & 1.0 & \bf{0.72} & 1.17 & 2.00 & 3.02 & \bf{3.89} \\
\hline
{\footnotesize $f_{11}(x) = (x + 1)^3 - 1$}  & \hfil 9.0 & \hfil 8.1 & \hfil \bf{6.3} & \bf{100} & \hfil 99 & \bf{100} & 1.0 & \bf{0.85} & 1.23 & 2.00 & 2.93 & \bf{4.08} \\
\hline
{\footnotesize $f_{12}(x) = (x-2)^3 - 10$}  & \hfil 9.3 & \hfil 8.9 & \hfil \bf{7.0} & \bf{100} & \hfil 98 & \bf{100} & 1.0 & \bf{0.96} & 1.34 & 2.00 & 2.99 & \bf{3.92} \\
\hline
{\footnotesize $f_{13}(x) = (x + 5/4) e^{(x + 5/4)^2} - \sin(x + 5/4)^2 + 3\cos(x + 5/4) + 5$}  & 11.6 & \hfil 9.4 & \hfil \bf{8.5} & \hfil 96 & \hfil 80 & \hfil \bf{97} & 1.0 & \bf{0.80} & 1.34 & 2.00 & 3.02 & \bf{3.64} \\
\hline
{\footnotesize $f_{14}(x) = x + \sin(2/x)  x^2$}  & 10.9 & 10.3 & \hfil \bf{8.0} & \hfil \bf{99} & \hfil 74 & \hfil 81 & \bf{1.0} & 1.08 & 1.28 & 2.00 & 2.99 & \bf{3.99} \\
\hline
\end{tabular}
\caption{\label{table:2}
Annealing effect on iterative root-finding. The numerical performances of fixed $\beta$ versus annealing $\beta_n$ are compared. The evaluations include: (i) the average number of iterations required for convergence (Iter. number); (ii) the percentage of initial points that converge (Conv. percentage); (iii) the relative computation time compared to the $\beta=0$ case (Comp. time); and (iv) the order of convergence (Conv. order). Note that the fixed $\beta=0$ and $\beta=1$ represent the original and the modified Newton-Raphson method, respectively. The iteration-dependent $\beta_n$ is formulated in Eq.~(\ref{eq:possible}). The convergence order is estimated by Eq.~(\ref{eq:ACOC}). In each category, the best performance—shortest iteration number, highest convergence percentage, least computation time, and largest convergence order among different $\beta$ values—is highlighted in bold.}
\label{tab:annealing}
\end{table*}

\subsection{Annealing} \label{anneal}

Indeed, the physical intuition gained from considering Eq.~(\ref{eq:hightemp}) as a high-temperature expansion suggests that annealing the temperature during iterations could improve the performance of the root-finding algorithm.
In our physical analogy, the auxiliary parameter $\beta$ in Eq.~(\ref{eq:iteration}) can be interpreted as the inverse of temperature in statistical mechanics.
With the high-temperature analogy in mind, it is natural to view the original Newton-Raphson update as the high-temperature limit ($\beta=0$), where only the linear Taylor expansion of the function is used, making the hopping term dominant.
In contrast, the Adomian update represents the low-temperature limit ($\beta=1$), where the exact minimum encoded in the $\beta$-dependent term is balanced with the linear hopping term.
As the temperature is decreased, the importance of the ultralocal exact root should increase. While exploration over as broad an area as possible is essential at the beginning, focusing on the root location by gradually decreasing the temperature should help in getting to the actual root. In fact, if we increase the inverse temperature from $\beta = 0,$ which corresponds to Newton-Raphson, to a final value $\beta$ after even one step, we find that the total number of iterations needed decreases and the dependence of the number of iterations on the final value of $\beta$ becomes much smoother as well. However, if $\beta$ is ramped up too gradually, there is no benefit to be had for most of the test functions we considered.
This annealing approach is similar to multipoint methods in that the update formulation changes with each iteration step~\citep{petkovic2010families}.

Let us now derive an annealing schedule from first principles. 
First, we introduce shorthand notation:
\begin{align*}
    f_{n} &\equiv f(x_{n}), \\ 
    f'_n &\equiv f'(x_n), \\
    \hat{f}_{n+1} &\equiv f(\hat{x}_{n+1}), \\
    \hat{f}'_{n+1} &\equiv f'(\hat{x}_{n+1}). 
\end{align*}
Then, our new update is represented as
\begin{align}
    \hat{x}_{n+1} &= x_{n} - \frac{f_{n}}{f'_{n}}, \nonumber \\
    x_{n+1} &= \hat{x}_{n+1} - \beta \frac{\hat{f}_{n+1}}{f'_{n}},
\end{align}
where $\hat{x}_{n+1}$ is the Newton-Raphson update.
A Taylor expansion of $f(x_{n+1})$ with the second equation gives the following relation:
\begin{align}
    f(x_{n+1}) &= f \left(\hat{x}_{n+1} - \beta \frac{\hat{f}_{n+1}}{f'_{n}} \right) \nonumber \\
    &\approx f(\hat{x}_{n+1}) -f'(\hat{x}_{n+1}) \beta \frac{\hat{f}_{n+1}}{f'_{n}}.
\end{align}
This implies 
\begin{equation}
    f^2_{n+1} =  \left[\hat{f}_{n+1} - \beta \frac{\hat{f}'_{n+1}\hat{f}_{n+1}}{f'_{n}}\right]^2
\end{equation}
with our shorthand notation.
Therefore,
\begin{equation}
\label{eq:inequality}
    f^2_{n+1} \le \hat{f}_{n+1}^2,\ \ \mathrm{if}\ \ \left|1- \beta \frac{\hat{f}'_{n+1}}{f'_{n}}\right| \le 1.
\end{equation}
In other words, our update makes $f^2(x_{n+1})$ smaller than the Newton-Raphson update $f^2(\hat{x}_{n+1})$, when $\beta \hat{f}'_{n+1}/f'_{n}$ is small enough to justify the Taylor expansion of ${f}$ around $\hat{x}_{n+1}.$ 
Now, suppose we take 
\begin{equation}
    \beta_n \equiv \frac{2{f_n'}^2}{{\hat{f}'^2_{n+1}} + {f_n'}^2} \ge 0.
    \label{eq:possible}
\end{equation}
This choice guarantees Eq.~(\ref{eq:inequality}) because
\begin{equation}
    1- \beta_n \frac{\hat{f}'_{n+1}}{f_n'} =  \frac{\left({{\hat{f}'_{n+1}} - {f_n'}}\right)^2}{{\hat{f}'^2_{n+1} + {f_n'}^2}} \le 1.
\end{equation}
Therefore, the annealing schedule in Eq.~(\ref{eq:possible}) suggests an approach to setting $\beta_n$ depending on the derivatives of the function at the previous value of $x_n$ and the Newton-Raphson update value $\hat{x}_{n+1}.$ 
Equation~(\ref{eq:possible}) has appropriate limits of $1$ or $0$ when $f'_n$ is large or small depending on the value of ${\hat{f}_{n+1}}'$ of course. When $f_n'$ is large and $x_n$ is proximate to roots, we expect the two derivative values are close, so $\beta_n \approx 1$ with $|1- \beta_n {\hat{f}'_{n+1}}/{f_n'}| \le 1$.
When it is small at the singularity, however, $\beta_n \approx 0$ with a finite derivative value $\hat{f}'_{n+1}$.

Table \ref{tab:annealing} summarizes the results for fixed values of $\beta$, and for the annealing schedule $\beta_n$ defined in Eq.(\ref{eq:possible}). The annealing schedule enhances the effectiveness of the $\beta$ value near roots. The estimated order of convergence near a root is $q_{n} \simeq 4$ for most functions. 

The mean number of iterations has been computed for a $1000\times1000$ regularly spaced grid using the same methods and parameters as in Sec.~\ref{sec:results}. The best results have been highlighted in bold. There is no function for which the Newton-Raphson algorithm performs better than the annealing schedule. We caution that the derived annealing schedule in Eq.~(\ref{eq:possible}) does require an additional derivative evaluation $f'(\hat{x}_{n+1})$. This extra computational cost is reflected in the total time per iterations. The time per point is indeed higher despite the reduction in the number of iterations due to the added derivative. A possible optimization is to approximate this derivative with first order approximations using the previously evaluated functions. 

In the supplemental Julia scripts available at \texttt{https://github.com/awage/RootFinding}, computations are provided for all the standard test functions listed in Tab.~\ref{table:1}.

\section{Discussion}\label{end}
Root-finding is a technique that is central to many quantitative science and engineering problems. This paper has investigated a new root-finding approach that has the potential to offer several benefits: improved efficiency and enhanced numerical stability, broader applicability to complex nonlinear functions that abound in real-world problems, and more thorough investigation of possible roots because of higher basin entropy, all while not increasing the order of derivatives needed for application. More fundamentally, the physical picture of the variation in the root basins of attraction as the temperature is decreased and the cohomological field theory explication of roots as transcendental cohomology classes are new contributions to the theoretical underpinnings of numerical analysis and algorithms. We have just scratched the surface of the applications of this technique, as is evident from the connection we made to the Adomian method, which has been used for solving nonlinear, and even stochastic, ordinary and partial differential equations. This connection implies that our approach could be applied to all these problem areas as well.

In particular, the homotopy analysis method~(HAM) is a well-known approach to solving nonlinear problems due to \citet{liao2003beyond}, motivated by topological homotopy theory. HAM unifies the Adomian method and a host of other numerical methods using an auxiliary parameter that constructs a homotopy to handle the nonlinear nature of the problem. This so-called convergence control parameter is then used to show convergence of a series solution. As HAM can be combined with spectral or Pad\'e approximation methods, it would be interesting to see if HAM could be combined with our approach, Eq.~(\ref{eq:beta}), to completely avoid a series expansion.

\section*{Acknowledgement}
V.P. thanks Christopher Kim for helpful discussions. This work was supported by the Intramural Research Program of the National Institutes of Health, NIDDK (V.P.), and by the Creative-Pioneering Researchers Program through Seoul National University, and the National Research Foundation of Korea (NRF) grant (Grant No. 2022R1A2C1006871) (J.J.)

\appendix
\section{Topological formulation of root-finding}
\label{topological}
We formulate a rigorous topological formulation of root-finding to show that root-finding iterations are explicit representatives of the same cohomology class. The derivation below is self-contained but we note that it is an extension of cohomological quantum field theory~\citep{Witten:1990bs}. 

We introduce two anti-commuting variables, $b$ and $c,$ with the following properties:
\begin{equation}
    bc = - cb,\ b^2 = 0,\ c^2 = 0,
\end{equation}
and a Grassmann integration over these variables defined by
\begin{align}
    \int\dd b 1 = 0,\ \int \dd b b = 1, \int\dd c 1 = 0,\nonumber \\ \int \dd c c = 1, \int\dd b c = 0,\ \int \dd c b =0 , 
\end{align}
extended by linearity and all the standard properties of integration.
We also introduce an anti-commuting derivation, $s,$ (called the BRST operator in the physics literature) as follows:
\begin{equation}
    sx = c, \ sc = 0,\ sb = \lambda,\ s\lambda = 0 
\end{equation}
where $\lambda$ is a commuting variable for consistency with the fact that combinations of anti-commuting quantities, in this case $s$ and $b,$ are commuting quantities. It is trivial to verify that $s^2 = 0,$ which implies that $s$ can be used to define a cohomology theory.

Now we define a gauge fermion depending on a parameter $g$ as follows: 
\begin{equation}
    \Psi = b(if(x) - g^2\lambda/2)
\end{equation}
where $f$ is the function for which we want to find roots.
The action of $s$ on $\Psi$ gives
\begin{equation}
    s\Psi = \lambda(if(x) - g^2\lambda/2) - ibf'(x)c
\end{equation}
where the negative sign of the second term arises when the derivation $s$ anti-commutes past $b.$

Consider the measure 
\begin{equation}
    \dd \mu \equiv \dd x\dd\lambda\dd b\dd c \exp(s\Psi) g/\sqrt{2\pi}.
\end{equation}
which is invariant under the action of $s.$ If we integrate a general function $h(x,c,b,\lambda)$ with respect to this measure, the value of the integral will depend on the parameter $g.$ However, suppose that $h$ satisfies $sh = 0.$ Then a change in $g$ will lead to a change in the value of the integral of the form 
\begin{equation}
   \delta {\cal I} = \int \dd\mu (- sc\lambda/2) h
\end{equation}
but the invariance of the measure under the action of $s$ allows an integration by parts to give us
\begin{equation}
    \delta {\cal I} = \int \dd\mu  (c\lambda/2) sh = 0.
\end{equation}
Thus, for very specific functions that satisfy $sh = 0,$ the value of their integral does not depend on the parameter $g.$ Moreover, if we change $h$ to $h + sk,$ where $k$ is an arbitrary function of the variables, we see by the same argument that the value of the integral of $h$ does not depend on $k.$ Functions that satisfy $sh = 0$ are called closed functions and functions that are of the form $sk$ are called exact functions.

These are the axioms that define a cohomology theory, and therefore closed functions fall into equivalence classes, defined up to the arbitrary addition of exact functions.

How is this formal development relevant for finding roots? A formal Laplace approximation to the measure gives the `equations of motion'
\begin{equation}
    x_\root:f(x_\root) = 0,\ \lambda=0,\ cf'(x_\root) =0,\ f'(x_\root)b = 0.
\end{equation}
Thus, this measure is localizing at the roots of $f.$
Integrating over $\lambda$ in
\begin{equation}
   {\cal Z} = \int \dd\mu \exp(s\Psi)
\end{equation}
gives 
\begin{equation}
    \hat {\cal Z} = \int \dd x\dd b\dd c\exp(- f(x)^2/2g^2 -icf'(x)b)
\end{equation}
where we should now expect to use the $\lambda$ equation of motion to verify that $s^2 = 0.$ The striking thing about this expression is that if we make $g$ small, we see that any closed function will depend only on the value of $x$ at which $f(x) = 0.$ Therefore, this formulation of our physical approach exhibits root-finding as a mathematically rigorous cohomology theory. Any nontrivial cohomology class in this problem will have to be a transcendental function, since it is easily demonstrated that the cohomology class of any polynomial function of the variables is guaranteed to be trivial.

Let us now demonstrate that the cohomology operator $s$ can help us check if an expression for a root of $f(x)$ is, or is not, $g$ independent. Consider the formal series 
\begin{equation}
    x_\root \equiv x - h(x) - h(x - h(x)) - \ldots
\end{equation}
where we have defined $h \equiv f/f',$ and this is simply the series of corrections of the Newton-Raphson algorithm.
We claim that $sx_\root = 0.$ This is explicitly verified as follows:
\begin{equation}
    sx_\root = c - \frac{f'(x)}{f'(x)}c + \frac{f(x)}{f'(x)^2}c + ...
\end{equation}
where the first two terms obviously cancel. This general pattern continues because $s$ acting on the $n^{\mathrm {th}}$ term gives
\begin{equation}
    s(x_n - x_{n-1}) = \frac{f(x_{n-1})}{f'(x_{n-1})^2} sx_{n-1} - sx_{n-1}.
\end{equation}
As the previous $n-1$ terms sum up precisely to $sx_{n-1}$ because this is a telescoping series, we have a cancellation. Since the equation of motion is $f(x_\root) = 0,$ we see that if the series converges, it is a representative of the cohomology class defined by the root. 

In particular, some algebra shows that our $\beta$ dependent variations on the representative series also have the same type of telescoping structure: 
\begin{widetext}
\begin{equation}
    s(x_n - x_{n-1}) = \frac{f(x_{n-1})}{f'(x_{n-1})^2} \left[1 + 
    \beta \left(\frac{{f(\hat x_{n-1})}}{f(x_{n-1})} - \frac{{f'(\hat x_{n-1})}}{f'(x_{n-1})} \right)\right]sx_{n-1} - sx_{n-1}.
\end{equation}
\end{widetext}
where we have defined $\hat x_{n-1} \equiv x_{n-1} - h(x_{n-1}).$ Thus, independent of $\beta,$ this defines a cohomology class. To see the promised simplification produced by our choice of coefficient in Eq.~(\ref{eq:xieom}), a little scrutiny reveals that if the series is converging, the term ${{f'(\hat x_{n-1})}/{f'(x_{n-1})}} \rightarrow 1,$ which leads to a cancellation at $\beta = 1, $ giving
\begin{equation}
    s(x_n - x_{n-1}) = \frac{f(\hat x_{n-1})}{f'(x_{n-1})^2} sx_{n-1} - sx_{n-1}.
\end{equation}
As $\hat x_{n-1}$ is the Newton-Raphson iterate of $x_{n-1}$, this suggests faster convergence, provided of course that the assumptions made in this simplification are valid.
This cohomological explanation also removes any mystery in the invariance of the roots no matter what value of $g$ we use for approximating the integral: we simply have to ensure that the expression for the root is a closed function for the cohomology operator $s.$ We show in section~\ref{faster} that the faster convergence suggested above can be explicitly demonstrated  with the usual derivation used to show the quadratic convergence of the Newton-Raphson algorithm.

\bibliographystyle{apsrev4-1}
\bibliography{newton}

\end{document}